\documentclass[twoside,letterpaper,draft,10pt]{amsart}

\usepackage{amsmath,amsthm,latexsym,xspace,amscd,amssymb, textcomp}

\usepackage[margin=1.8in]{geometry}

\title[Nil-Generated Algebras and Group Algebras]{Nil-Generated Algebras and Group Algebras\\
	 whose units satisfy a Laurent polynomial identity}

\dedicatory{ Claudenir Freire Rodrigues	\\
	         \vspace{0,2cm} 
            Departamento de Matemática, Universidade Federal do Amazonas, Amazonas,\\ Manaus, Brasil\\
            E-mail: claudiomao1@gmail.com }

\newtheorem{teo1}{Theorem}[section]

\newtheorem{lem1}[teo1]{Lemma}

\theoremstyle{remark}\newtheorem{remark}[teo1]{Remark}
\theoremstyle{remark}

\begin{document}

\begin{abstract}
Let $A$ be an algebra whose group of units $U(A)$ satisfies a Laurent polynomial identity (LPI). We establish conditions on these polynomials in such a way that nil-generated algebras and group algebras with torsion groups over infinite fields in characteristic $p>0$ have nonmatrix identities. We also determine, in the context of group algebras with arbitrary LPI for the group of units, the existence of polynomial identities.  
\end{abstract}

\maketitle

\section{\textbf{Introduction}}
	Let $N_A$ be the set of nilpotents of an $R$ algebra $A$. If $A$ is generated by $\{1\}\cup N_A $, we say that $A$ is nil-generated.
	
A Laurent polynomial  in the noncommutative variables $x_i, \, i=1,\ldots,l$ is a non-zero element in the group algebra  $RF_l$ over a ring  $R$  with free group $F_l = < x_1,...,x_l >$.    A $ P $ in $RF_l $  is said to be a Laurent Polynomial Identity  (LPI) for $U(A)$ (the group of units in A) if $P(c_1, ..., c_l)=0 $ for all sequence $c_1,..., c_l$ in $U(A)$. The special case $P = 1-w$ is called a group identity (GI) . The sum of the exponents of a word  $w=w(x_1,\ldots,x_n) \in F_l  $ in the variables $x_i$ is given by $\sum exp _{x_i}\, w$. Replacing $x_ i$ by $x^{-i}yx^i$ we can suppose the LPI is in two variables. A polynomial identity not satisfied by $M_2(R)$ is a nonmatrix identity. 

Polynomial identities and group identities have been studied with a focus on the conjecture of Hartley \cite{SKS}, which says that if 
\ $G$ is a torsion group, $K$ a field, and the group $U(KG)$  satisfies a group identity, then $KG$ satisfies a polynomial identity (PI). This was proved first over infinite fields by \cite{GSV}, established definitely by \cite{LIU}, and characterizations were given by \cite{OJA}  with LPI  instead of GI  to guarantee the existence of polynomial identities.  In addition, results were obtained for nil-generated unitary algebras in \cite{BRT}, confirming that the existence of GI for $U(A)$ is equivalent to the existence of a nonmatrix identity for $A$. In the first section of the present paper we treat these questions with  LPI for $U(A)$, obeying a condition that guarantees  the existence of a nonmatrix identity  for nil-generated unitary algebras  (Theorem \ref{n}), which is equivalent to the existence of a group identity for the group of units of these algebras. This condition will now be justified. For an LPI $P=\, a_1 + a_2w_2 +\cdots+a_n w_n $, with $a_1$ non-zero, the words $w_{i\geq 2}$ are the nonconstant words. According to  Amitsur--Levitzki \cite[Theorema V. 1.9]{P},  $U(M_n(K) \, n>1$, satisfies the LPI 
$ f(x_1,...,x_{2n})= S_{2n}\cdot(x_1\cdots x_{2n})^{-1}  $  where $S_{2n}$ is the  standard identity and $f$ is an  LPI with all non-constant words having sum of exponents zero at every variable. So if $K$ is infinite, then \cite[Lemm3.3]{LIU} implies that $U(M_n(K)$ can not satisfy a group identity. For this reason, we impose the requirement that the sums of the exponents of non-constant words in the LPI be non-zero in at least one of the variables. In the second section, this assumption is necessary for the determination of PIs for group algebras of a finite group over an infinite field (Theorem \ref{n2}.(2)). From the context of the proof of Theorem \ref{n2}.(1), we get a confirmation of the Hartley Conjecture with arbitrary LPI for the group of units $U(KG)$ with $G$ a torsion group and $K$ a field (Theorem \ref{k}). And finally, in the third section, we show the existence of a nonmatrix identity for a group algebra of a  torsion group over an infinite field in characteristic $p>0$ (Theorem \ref{y3}).    
 
\section{The Nil-Generated Case}

We need some results that adapt known classical ideas. The next two results are taken from \cite[Theorems 1.3--1.4]{BRT}.

\begin{teo1} \label{y1}
	Let $A$ be a nil-generated unitary algebra over a field of characteristic 0 and $U(A)$ its group of units. Then the following statements are equivalent:
	\begin{enumerate}
		\item[1.] $U(A)$ satisfies a group identity;
		\item[2.] $A$ satisfies a nonmatrix identity;
		\item[3.] $A$ is Lie soluble; 
		\item[4.] $U(A)$ is a soluble group.
	\end{enumerate}
	
\end{teo1}	 

\begin{teo1}\label{y2}
	Let $A$ be a nil-generated unitary algebra over a field of characteristic $p>0$ and $U(A)$ its group of units. Then the following statements are equivalent:
	\begin{enumerate}
		\item[1.] $U(A)$ satisfies a group identity;
		\item[2.] $A$ satisfies a nonmatrix identity;
		\item[3.] $A$ satisfies a polynomial identity $([x_1,x_2]x_3)^{p^t}=0$ for some $t$; 
		\item[4.] $U(A)$ satisfies the group identity $(y_1, y_2)^{p^t}=1$ for some $t$.
	\end{enumerate}
\end{teo1}

	From \cite[Theorem 2.1]{C} we have

\begin{teo1}\label{s1}
	
	Let $A$ be an algebra whose  group of units $U(A)$ admits an  LPI $P$  over a ring $R$ with unit  and such that the non-constant words $w_i$ in $P$ have sum of exponents  non-zero in at least one of the variables. Then there exists a polynomial $f \, \in \, R[X]$ with  degree $d$ determined by $ l=min \, \{\sum exp \,w_i\}$ and $r=max\{\sum exp \, w_i \}$  such that for all $ \, a,b,c, u $ in $A$ with $a^{2}=bc=0$,  $f(bacu)=0$. Thus, $bacA$  has a polynomial identity.
	
\end{teo1}

 From now one, unless otherwise stated, the LPI are as in Theorem \ref{s1}. And 
 according to \cite[Corollary 3.3]{C}, we have the following lemma.
 
 \begin{lem1}\label{s2}
 	
 	Let $A$ be a semiprime algebra over an infinite  commutative domain $R$. If $U(A)$ satisfies an LPI, then every idempotent element of $R^{-1} A$ is central and for all $b,\,c \in A$ such that $bc=0$, we have $bac=0$ for every $a$ nilpotent in A.			
 	
 \end{lem1}

A known result in the literature is
\begin{lem1}\label{r1}
	Suppose $f:R\to S$ is an epimorphism of rings with kernel a nil ideal. Then $f$ induces an epimorphism $ U(R)\to U(S)$.
\end{lem1}

Now we are ready to begin our adaptations.

 \begin{teo1}\label{s3}
	Let $A$ be an algebra over $R$ whose  group of units $U(A)$ admits an  LPI $P$. Then 
	
	\begin{enumerate}
		
		\item[1.] There are identities in $R[x]$ determined by $P$ for the group of units $U(A)$ and the Jacobson radical $J(A)$. These  are nonmatrix identities if $R$ is an infinite field.
		
		\item[2.] There is a polynomial $g(x) \in R[x]$ determined by $P$ such that for all $a,\,b \in A$ with $a^2=b^2=0$, $g(ab)=0$. And, conversely, the ideal $bacA$ where $a^2=b c=0$ in $A$ has a polynomial identity.
		
		\item[3.] Suppose $A$ is a nil generated algebra and $R$ is an infinite domain. Then the set $N_A$ of nilpotents in $A$ is a locally nilpotent ideal and $A$ has a  polynomial identity. This identity is nonmatrix if $R$ is an infinite field. And every finite set of nilpotent elements generates a nilpotent subalgebra.
		
	\end{enumerate}	
	
\end{teo1}

\begin{proof}

 Let  $P\,=\, a_1 + a_2w_2 +...+ a_t w_t$ be the LPI with $w_i $ in the form $$w_i=x_1^{r_1}x_2^{s_2} \cdots x_1^{r_k}x_2^{s_k}$$ with $k\geq 1, r_i,\, s_i $ integers.  
  
   If $ \sum exp \, w_i = \sum exp _{x_1}\, w_i + \sum exp _{x_2}\, w_i = 0 $ then we  substitute $ x_1= x_1^k$ or $x_2= x_2^k$ with $k > 1$ big enough.
	We get a new LPI in the form  $P\,=\, a_1 + a_2{w'}_2 +...+ a_r {w'}_r,$ where   
	$$ \sum exp \, {w'}_i = \sum exp _{x_1}\, {w'}_i + \sum exp _{x_2}\, {w'}_i = k \sum exp _{x_1}\, w_i + \sum exp _{x_2}\, w_i $$  or $$\sum exp \, {w'}_i = \sum exp _{x_1}\, w_i +  k\sum exp _{x_2}\, w_i. $$     In any  case that is not zero. So we suppose $ \sum exp \, w_i $ is not zero for all $i = 2,...,t$.		 
	
	Hence, $$P(\alpha, \alpha)\,=\, a_1 + a_2 \alpha ^l +...+a_i \alpha ^{\sum exp \, w_i} +...+ a_t \alpha ^r =  0   $$ with all powers of $\alpha$ not zero for all $ \alpha \in U(A)$ and  $l$ and $r$ are respectively the minimum and maximum of the sum of the exponents of the  $w_i's$.
		
	After collecting like terms, we will have a polynomial expression of the form	
	$$  f_0(\alpha) = a_1  +  b_2 \alpha ^l +...+ b_j \alpha ^r=  0 $$ where $ b_j's$ are the partial sums of $a_i's, \, i=2,...,t $.	
	In  general,  $P(1,1)\,=\, a_1 + a_2 +...+ a_t =0$, so $ a_2 +...+ a_t \neq 0 $, because $a_1 \neq 0$. So we can not have all $b_i's =0$ because in this case $a_2 +...+a_r = b_2+...+b_j=0$. It follows that $f_0$		
	is a non-zero polynomial over $R$ with integer exponents not all   necessarily positive. 
		
	If $l < 0$, we have that  $\alpha ^{-l}f_0(\alpha)= f_1(\alpha)=0$, where $f_1\neq 0$ is a polynomial identity for $U(A)$  over $R$ of degree at most $r-l$. Suppose $l>0$. Then $f_0$ is a polynomial identity for $U(A)$. Now, with $a$ in $J(A)$,  $1 \,+ \,a $ is a unit and $ f_0(1 \,+\,a) = 0 $. It follows that $g(x)=f_0(1 \,+\,x)$ is an identity for  $J(A)$. For all $\lambda \in R$ the minimal polynomial of the matrix $\lambda I_{2\times2}$ where $I_{2\times2}$ is the identity matrix is $x - \lambda$. If a polynomial in $F[x]$ is a matrix identity, then $\lambda$ is a root of it and $R$ would be a finite field. This finishes the first part.
	
	In the sequel, by Theorem \ref{s1}, for the case $b=c$, there is a nonzero polynomial  $f(x) \in R[x]$ with $f(babu)=$ for all $u\in A$. Now take $u=1$, $af(bab)=g(ab)=0$. The converse is known, see \cite[lemma 3.1]{OJA}. Let $g_1(x) \in R[x]$ be such that $g_1(ab)=0$ if $a^2=b^2=0$ in  $A$. For all $u\in A$, $(cub)^2=0$, and then $bg_1(acub)acu=g(bacu)=0$. 
	
	Finally, according to \cite{J} $A/L$, where $L$ is the lower nil radical of $A$ (which coincides with the intersection of all prime ideals of $A$), is a semiprime ring and $L$ is a locally nilpotent ideal in $A$. It follows that $ U(A) \to  U(A/L) $ is an epimorphis by Lemma \ref{r1} and so $U(A/L)$ satisfies $P$. Hence, using Lemma \ref{s2}, for all $a,\, x \in A/L$ with $a^2=0$ and $x$ nilpotent, we have $axa=0$. And by induction, $ax_1\ldots x_na=0$ with the $x_i$'s nilpotents in $A$. As $L$ is a nil ideal, $A/L$ is nil generated, so $a(A/L)a=0$ and   $(A/L)a(A/L)$ is a nilpotent ideal. From the fact that $A/L$ is semiprime, it follows that  $a=0$, which gives us that there are no nilpotents in $A/L$. Then $L$ is the set of nilpotents in $A$ and $A=R + L$. As for all $x,\,y \in A,\, [x,y] \in L \subset J(A)$, $A$ satisfies the polynomial identity $f([x,y])$ where $f$ is the identity determined by the LPI in the first part we have proved. With	 $a=\lambda e_{12},\,b=e_{21}$ in $M_2(R)$, we have that
	 $$f([a,b])= \left[\begin{array}{cc} 
		            f(\lambda)&0\\
		            0&*\\
	\end{array}\right]$$ for all $\lambda \in R$. It follows that $f([x,y])$ is a nonmatrix identity if $R$ is an infinite field.
	 With $L$ locally nilpotent we finish the proof.
		   
	\end{proof}

	\begin{teo1}\label{n}
		
		Let $A$ be a nil-generated unitary algebra over an infinite field $F$ with the group of units $U(A)$ satisfying an LPI. Then $A$ satisfies a nonmatrix identity. And the following are equivalent in characteristic $0$:
		       
    \begin{enumerate}
    	\item[1.] $U(A)$ satisfies a group identity;
    	\item[3.] $A$ is Lie soluble; 
    	\item[3.] $U(A)$ is a soluble group.
    \end{enumerate}

	And the following are equivalent in characteristic $p>0$:
	
	\begin{enumerate}
		\item[1.] $A$ satisfies a nonmatrix identity;
		\item[2.] $U(A)$ satisfies a group identity;
		\item[3.] $A$ satisfies the polynomial identity $([x_1,x_2]x_3)^{p^t}=0$ for some $t$; 
		\item[4.] $U(A)$ satisfies the group identity $(y_1, y_2)^{p^t}=1$ for some $t$.
	\end{enumerate}

\end{teo1}

\begin{proof}
	Immediate consequence of Theorems \ref{y1}, \ref{y2} and \ref{s3}.
\end{proof}
 
\section{Group Algebras and Standard Identities}

According to \cite{OJA}, an $R$-algebra $A$ has  the property $\mathrm{P_1}$ with respect to a nonzero polynomial $g[x]\in R[x]$ if $g(ab)=0$ for every $a,\,b \, \in R $, with $a^2=b^2=0$. And a Laurent polynomial $P$ over $R$ has the property $\mathrm{P}$ if every $R$-algebra for wich $L$ is an LPI of $U(A)$ has the property $\mathrm{P_1}$ with respect to some nonzero polynomial.

\begin{lem1}\label{Gon84}
	Let $D$ be a noncommutative division ring, finite dimensional over its center. Then $U(D)$ contains a free subgroup of rank two.
	
\end{lem1}	

\begin{proof}
	See \cite{Gon84}.
\end{proof}

 \begin{teo1}\label{n2}
 	Let $G$ be a finite group and $F$ a field of characteristic  $p\geq 0$ with $U(FG)$ satisfying  an LPI $P=\, a_1 + a_2w_2 +\cdots+a_n w_n $.   
 \end{teo1}

\begin{enumerate}
	 \item [1.] If $P$  is an arbitrary LPI, then there exist positive integers $m,\,t$ such that the group algebra  $FG$ satisfies the polynomial identity  $(S_{2m})^t=0$. And if $F$ is of charactheristc 0 or characteristic $p$ where $p\not |\, |G|$, then $FG$ satisfies the standard identity $S_{2m}$.
	 
	\item [2.]  If $F$ is infinite, then $FG$ satisfies the nonmatrix identity $ (S_2)^t=[x,y]^t=0 $.\\ When $charF=0$, $G$ is an abelian group.\\	
	If $charF=p>0$, $U(FG)$ satisfies the group identity $(x,y)^{p^k}=1$ for some $k$  and  in this case, $G$ is a $p$-abelian group. In particular, if $G$ is a $p'$-group, then $G$ is abelian.

\end{enumerate}

\begin{proof} 
	 As $G$ is finite, we have that $FG$  and $ FG/J$
		 are  Artinian and the Jacobson radical  $J=J(FG)$ is nilpotent. By Lemma \ref{r1} the epimorphism $FG\to  FG/J$ induces the epimorphis   $\pi: U(FG)	 \to  U(FG/J) $  and so $U(FG/J)$ satisfies  $P$.
		 
		  As $ FG/J$ is 
	 	 semisimple, by the Wedderburn--Artin Theorem \cite[pg 35]{Lam} $ FG/J \cong \oplus_{i=1}^r M_{n_i}(D_i) $ with $D_i$  division rings and the center $Z_i$ of $D_i$ contains $F$,  
    such that $[D_i:Z_i]\leq [D_i:F]<\infty$ because $G$ is finite. Assume that  $D_i$ is noncommutative. Then, according to Lemma \ref{Gon84}, $U(D_i)$  contains a non-abelian  free group $H$ of rank 2, $H=<g_1,g_2>$ where $\{g_1,g_2\}$ is a free generator set of $H$.  With the units $g_1$ and $g_2$, as $U(D_i)$ satisfies $P$, then we have
    
   $$P(g_1,g_2)= \, a_1 + a_2w_2(g_1,g_2) +\cdots+a_n w_n(g_1,g_2)=0 $$ 
   
   For every reduced word $w_i$, $w_i(g_1,\,g_2)$ is unique and not trivial in $H$. The second equality above tells us that  every $a_i$ is zero  because the $w_i(g_1,g_2)$ are linearly independent in the subalgebra $FH$. With this contradiction, $D_i$ is a field for each $i$. And with $m=max\{n_i \}$,  $ FG/J $ satisfies $S_{2m} $. It follows that for all $\alpha_1, \ldots,\alpha_{2m}$ in $FG$, $S_{2m}(\alpha_1, \ldots,\alpha_{2m})$ is in $J$ and so there exists a $t>0$ with $(S_{2m}(\alpha_1, \ldots,\alpha_{2m}))^t=0$. And if the characteristic is 0 or $p>0$ with $p\not | \,|G|$, then by Maschke's Theorem, $FG$ is semisimple satisfying $S_{2m}$. 
   
       For the second part, suppose  that  every non-constant word $w_i$ has  it  sums of exponents non zero in at least one of the variables. As $FG/J$ is semiprime, by Lemma \ref{s2} the idempotents are central and thus $n_i=1$ for every $i$.
     It follows that $ FG/J \cong \oplus_{i=1}^r D_i $ and is commutative. Thus for all $x,y \in FG$, $[x,y] \in J$ and $[x,y]^t=0$. As we saw in the proof of Theorem \ref{s3}(3), over an infinite field this is a nonmatrix identity.

       Moreover for all $x,y \in U(FG)$
       $$\pi(x,y)=(\pi(x),\pi(y))=1=\pi(1), $$
        and consequently $$(x,y) =1 modJ. $$ 
        Thus $$G'\subseteq 1 + J, $$ and it follows that  $$\bigtriangleup (G') \in J $$ where $\bigtriangleup (G')$ is the augmentation ideal of $FG'$. And so $\bigtriangleup (G')$ is nilpotent. Hence, $G'$ is a finite $p$-group when $charF=p >0$ and it is then immediate that $G$ is abelian if $G$ is a $p'$-group. If $charF=0$, $\bigtriangleup (G')$ is not nilpotent and $G'=\{1\}$, that is, $G$ is abelian.    
     
     To finish the second part, with $charF=p>0$. As we saw before, $(x,y) =1 modJ $. So there is a $k$ such that $(1-(x,y))^{p^k}=0$ and $U(FG)$ satisfies $(x,y)^{p^k}=1.$ The conclusion follows by \cite[Lemma 2.3]{GSV}.

\end{proof}

 In the next results, there is no restriction on the LPI.
 
\begin{teo1}\label{k}
	Let $KG$ be the group algebra of the torsion group $G$ over the field $K$. If $U(KG)$ satisfies an LPI, then $KG$ satisfies a PI.
\end{teo1}

\begin{proof}
	
	Let $F=K[\alpha,\, \beta : \alpha^2=\beta^2=0]$ be the free algebra in two noncommutative variables $\alpha,\, \beta$ relative to the relations $\alpha^2=\beta^2=0$. According to \cite{OJA}, the group of units $U=U(F)$ contains a non-abelian  free group $H$, say us with $g_1,\ldots ,g_n$  free generators of $H$. Arguing as in the proof of Theorem \ref{n2}.(1), we conclude that  $U(F)$ can not satisfy $P$. It follows from  \cite[Theorem 1.1]{OJA} that $KG$ satisfies a polynomial identity.
	 
\end{proof}	

\begin{remark}\label{f}
		According  to Theorem \ref{s3}.(2), if $ P=\, a_1 + a_2w_2 +\cdots+a_n w_n $ is an LPI   with all non-constant words having sum of exponents non-zero in at least one of the   variables, then $P$ satisfies the property $\mathrm{P_1}$ with respect to a polynomial $g(x) \in R[x].$ But for this conclusion, this restriction in general is not necessary because with the last proof we are guaranteed to have more information. As we saw, no Laurent polynomial can be an LPI for the units $U$ of the free algebra $K[\alpha,\, \beta : \alpha^2=\beta^2=0]$ in the noncommutative variables $\alpha,\, \beta$ relative to the relations $\alpha^2=\beta^2=0$. With this, we will be able to remove the word ``if'' in \cite[Lemma 1.3]{OJA}.
\end{remark}
	
	\begin{lem1}\label{u}
		Let $f$ be a Laurent polynomial over a field $K$. Then $f$ is not an LPI of $U$. And there is a non-zero polynomial  $g \in K[x]$ such that every $K$-algebra $B$ for which $f$ is an LPI of $U(B)$ has the property $\mathrm{P_1}$ with respect to $g$.		
	\end{lem1}

\begin{proof}
	We know the first part by  Remark \ref{f} and the equivalency is by \cite[Lema 1.3]{OJA}.
\end{proof}

 \begin{teo1}
	 	
	 	If $K$ is a field of characteristic $p$ and $G$ is a locally finite $p'$-group with $U(KG)$ satisfying a LPI, then $KG$ satisfies a standard polynomial identity. 
	 
 \end{teo1}		

\begin{proof}
	
	This is by Lemma \ref{u} and \cite[Proposition 1.4]{OJA}.
	
\end{proof}

Now we return to imposing the restrictions about the exponents on the  LPI.
  
  Nil rings give rise to multiplicative groups, called adjoint groups. If $R$ is nil, then its adjoint group $R^0$ consists of the elements
  of $R$ as the underlying set, with multiplication given by $r \circ s = r+s+rs$. In the case of a nil-algebra $A$ over a field $F$, it is known that we can embed $A$ into the unitary nil-generated algebra $A_1= F + A $ and by mapping $ a+r \to (a, a^{-1}r)$ where $a$ is in $F \backslash  \{0\} $ and $r$ is in $A$ we have that $U(A_1) \cong (F \backslash  \{0\}) \times A^0$. 
    
  \begin{lem1}\label{s4}
  	Let $A$ be a nil-algebra over an infinite field $F$ with the adjoint group $A^0$ satisfying an LPI. Then $A$ satisfies a nonmatrix identity. Furthermore, $A$ satisfies a group identity.
  
  \end{lem1}

\begin{proof}
		As we saw, $U(A_1) \cong (F \backslash  \{0\}) \times A^0$. By Theorem \ref{s3}, $A_1$ satisfies a nonmatrix identity. The result follows by \cite[Corollary 1.5]{BRT}.
	
\end{proof}

	With this last adaptation, we are going to finish these notes with a version of \cite[Theorem 4.1]{BRT} in the context of LPIs.
  
 \begin{teo1}\label{y3}
 	
 	Let $FG$ be the group algebra of a torsion group over an infinite field $F$ of characteristic $p > 0$ whose group of units satisfies an LPI $P$. Then $FG$ satisfies a nonmatrix identity. And the following are equivalent.
 	\begin{enumerate}
 		\item [1.] $U(FG)$ satisfies a group identity;
  		\item [2.] $G$ contains a normal subgroup $A$ such that $G/A$ and $A'$ are finite and $G'$ is a $p$-group of finite exponent;
 		\item [3.] $[FG,FG]FG$ is nil of bounded index; and
 		\item [4.] $(U(FG))'$ is a $p$-group of finite exponent.
 	\end{enumerate}
 \end{teo1}
        
	\begin{proof}
		 
		 According to  \cite[Theorem 4.1]{BRT}, it is guaranteed that the result holds if $FG$ has a nonmatrix identity. With the subroup $P_G$ generated by the set of $p$-elements in $G$, $FP_G$ is nil-generated, $U(FP_G)$ satisfies $P$, and according to Theorem \ref{s3} the set of nilpotent elements $N_{FP_G}$ forms a locally nilpotent maximal ideal in $FP_G$. Thus $FP$ coincides with the augmentation ideal of $FP_G$ and $P_G$ is a locally finite $p$-group. As $P_G$ is normal in $G$, $N(FP_G)FG$ is a locally nilpotent ideal in $FG$. By the mapping $a \to 1+a$ we  identified the adjoint group of $N(FP_G)FG$ with a subgroup of $U(FG)$. It follows from Lemma \ref{s4} that $N(FP_G)FG$ has a nonmatrix identity.
		 
		  And as $N_{FP_G}$ is the kernel of the projection $FG \to F(G/P_G)$ then, by Lemma \ref{r1}, $F(G/P_G)$ satisfies $P$. Now we are going to show that $G/P_G$ is abelian. With this, we will finish the proof because $FG$ will be a commutative extension of an algebra satisfying a nonmatrix identity. 
		 
		 As $\overline{G}=G/P_G$ is a $p'$-group, then, by \cite[p. 130--131]{P}, $FP$ is semiprime and according to Lemma \ref{s2}, for all $g \in \overline{G}$ the idempotent $n^{-1}$ $(1 + g + g^2 +\ldots +g^{n-1})$ is  central in $F\overline{G}$. Hence $\langle g \rangle$ is a normal cyclic group  in $\overline{G}$. Thus $\overline{G}$ is an abelian or Hamiltonian group. As $\overline{G}$ is Hamiltonian, it contains a copy of $Q_8$, the quarternion group of order 8, and $FQ_8$ satisfies $P$. By Theorem \ref{n2}.(2), $Q_8$ would be abelian. So, we have $\overline{G}$ abelian. This completes the proof..
			
	\end{proof}	

\printindex


\begin{thebibliography}{99}
	
	\bibitem{GJV94} A. Giambruno, E. Jespers, and A. Valenti, Group identities on units of rings, \textit{Arch. Math}. (Basel)	63(4), 291--296 (1994). 
	
	\bibitem{GSV} A. Giambruno, S. K. Sehgal, and A. Valenti, Group algebras whose units satisfy a group identity, \textit{Proc. Amer. Math. Soc}. 125, 629--634 (1997). 
	
	\bibitem{OJA} Broche, O., Gonçalves, J. Z. and del Río, Á. Group algebras whose units satisfy a Laurent polynomial identity. \textit{Arch. Math}. 111, 353–367 (2018). 
	
	\bibitem{C} C. F. Rodrigues. A note in algebras with Laurent polynomial identity. arXiv.2012.01500v1 [math. RA].
	
	\bibitem{P} D. S. Passman, \emph{The Algebraic Structure of Group Rings}, Dover, New York, 2011.
	
	\bibitem{P1} D. S. Passman, Group algebras whose units satisfy a group identity II, \textit{Proc. Amer. Math.  Soc}. 125(3), 657--662 (1997).
		
	\bibitem{J} N. Jacobson, \emph{PI-Algebras: An Introductiong}, Springer Lecture Notes in Math., Vol. 441, Springer-Verlag, Berlin, 1975.
	
	\bibitem{Gon84} J. Z. Gonçalves, Free subgroups of units in group rings, \textit{Canad. Math. Bull}. 27(3), 309--312 (1984). 
		
	\bibitem{LIU}L. Chia-Hsin, Group Algebras With Units Satisfying a Group Identity, \textit{Proc. Amer. Math. Soc}. 127, 327--336 (1999).
	
	\bibitem{SKS} S. K. Sehgal \emph{ Units in Integral Group Rings}, Essex, 1993.	
	
	\bibitem{Lam} T. Y. Lam \emph{ A First Course in Noncommutative Rings}, Springer-Verlag, New York, 1991.
	
	\bibitem{BRT} Y. Billig, D. Riley, and V. Tasic Nonmatrix varienties and nil-generated algebras whose units satisfy a group identity, \textit{J. Algebra} 190, 241--252 (1997).
\end{thebibliography}
\end{document}